\documentclass[12pt]{article}
\usepackage{amssymb,amsfonts,amsmath, psfrag,eepic,colordvi}
\usepackage[enableskew]{youngtab}

\topskip  
\makeatletter
\newcommand{\rmnum}[1]{\romannumeral #1}
\newcommand{\Rmnum}[1]{\expandafter\@slowromancap\romannumeral #1@}
\makeatother

\numberwithin{equation}{section}
\newtheorem{theorem}{Theorem}[section]

\newtheorem{lemma}[theorem]{Lemma}

\newtheorem{example}[theorem]{Example}

\newtheorem{observation}[theorem]{Observation}
\begin{document}
\parskip 6pt

\pagenumbering{arabic}
\def\sof{\hfill\rule{2mm}{2mm}}
\def\ls{\leq}
\def\gs{\geq}
\def\SS{\mathcal S}
\def\qq{{\bold q}}
\def\MM{\mathcal M}
\def\TT{\mathcal T}
\def\EE{\mathcal E}
\def\lsp{\mbox{lsp}}
\def\rsp{\mbox{rsp}}
\def\pf{\noindent {\it Proof.} }
\def\mp{\mbox{pyramid}}
\def\mb{\mbox{block}}
\def\mc{\mbox{cross}}
\def\qed{\hfill \rule{4pt}{7pt}}
\def\block{\hfill \rule{5pt}{5pt}}

\begin{center}
{\Large\bf     Alternating
permutations with restrictions  and standard Young tableaux } \vskip 6mm
\end{center}

\begin{center}
{\small    Sherry H.F. Yan  and Yuexiao Xu\\[2mm]
 Department of Mathematics, Zhejiang Normal University, Jinhua
321004, P.R. China
\\[2mm]
\rm huifangyan@hotmail.com
  \\[0pt]
}
\end{center}

\noindent {\bf Abstract.}   In this paper, we   give
 bijections between the set of $4123$-avoiding down-up alternating permutations
of length $2n$ and the set of standard Young tableaux of shape
$(n,n,n)$, and  between the set of $4123$-avoiding down-up alternating
permutations of length $2n-1$ and the set of shifted  standard Young
tableaux of shape $(n+1, n, n-1)$ via an intermediate structure of
Yamanouchi words.  Moreover, we get the enumeration of
$4123$-avoiding  up-down alternating permutations of even and odd
length by presenting  bijections between  $4123$-avoiding  up-down
alternating permutations  and  standard  Young tableaux.

\noindent {\sc Key words}: alternating permutation, standard Young tableau, shifted standard Young tableau.\\
\noindent {\sc AMS Mathematical Subject Classifications}: 05A05, 05C30.


\section{Introduction}
A permutation $\pi=\pi_1\pi_2\ldots\pi_n$ of length $n$ on
$[n]=\{1,2,\ldots, n\}$ is said to be an {\em up-down  alternating}
permutations if $\pi_1<\pi_2>\pi_3<\pi_4>\cdots$. Similarly, $\pi$
is said to be a {\em down-up} alternating permutation if
$\pi_1>\pi_2<\pi_3>\pi_4<\cdots$.  We denote by $\mathcal{UD}_n$ and
$\mathcal{DU}_n$ the set of up-down and down-up alternating
permutations of length $n$, respectively. Note that the {\em complement} map
$\pi=\pi_1\pi_2\ldots \pi_n \longmapsto (n+1-\pi_1)(n+1-\pi_2)\ldots
(n+1-\pi_n)$ is a bijection between the set $\mathcal{UD}_n$ and the
set $\mathcal{DU}_n$. Denote by $\mathcal{S}_n$  the set of all
permutations on $[n]$.  Given a permutation $\pi=\pi_1\pi_2\ldots
\pi_n\in \mathcal{S}_n$ and a permutation $\tau=\tau_1\tau_2\ldots
\tau_k\in \mathcal{S}_k$, we say that $\pi$ contains {\em pattern}
$\tau$ if there exists a subsequence $\pi_{i_1}\pi_{i_2}\ldots
\pi_{i_k}$ of $\pi$ that is order-isomorphic to $\tau$. Otherwise,
$\pi$ is said to  {\em avoid } the pattern $\tau$ or be {\em
$\tau$-avoiding}.

 The classic problem of enumerating permutations
avoiding a given pattern has received a great deal of attention and
has led to interesting variations. For a thorough summary of the
current status of research,
 see B\'{o}na's book \cite{bona}.
  As an interesting variation, Mansour \cite{mansour} studied alternating permutations avoiding a given pattern.
  Alternating permutations have the intriguing property \cite{mansour, stanley2} that for any pattern of length $3$,
  the number of alternating permutations avoiding that pattern is given by Catalan numbers. This property is shared by
  the ordinary permutations. This coincidence suggests that pattern avoidance in alternating permutations and in ordinary
  permutations may be closely related, which motivates the pattern avoidance in alternating permutations. Guibert and Linusson
   \cite{guibert} showed that doubly alternating Baxter permutations are counted by the Catalan numbers and Ouchterlony \cite{ouch}
   studied the problem of enumerating doubly alternating permutations avoiding patterns of length $3$ and $4$.
   Recently,  Lewis \cite{Lewis1} initiated the study of enumerating alternating permutations avoiding a given pattern of length
   $4$.   Let $\mathcal{UD}_n(\tau)$ and
$\mathcal{DU}_n(\tau)$ be the set of $\tau$-avoiding up-down and
down-up alternating  permutations of length $n$, respectively.
  Lewis \cite{Lewis1} provided bijections between the set $\mathcal{UD}_{2n}(1234)$ and standard Young tableaux of shape $(n,n,n)$, and
    between the set $\mathcal{UD}_{2n+1}(1234)$ and standard Young tableaux of shape $(n+1, n,
    n-1)$. By applying  hook length formula for  standard Young tableaux
    \cite{stanley1}, the number of $1234$-avoiding up-down
    alternating permutations of length $2n$ is given by ${2(3n)!\over
    n!(n+1)!(n+2)!}$, and the number of $1234$-avoiding up-down
    alternating permutations of length $2n+1$ is given by ${16(3n)!\over (n-1)!(n+1)!(n+3)!}
    $.      Using    the method of generating trees,
     Lewis \cite{Lewis2} constructed   recursive bijections between the set $\mathcal{UD}_{2n}(2143)$ and the set of standard Young tableaux
     of shape $(n,n,n)$, and between the set $\mathcal{UD}_{2n+1}(2143)$ and the set of shifted standard Young tableaux of shape $(n+2, n+1, n)$.
     Using computer simulations,  Lewis \cite{Lewis2} came up with several conjectures    that indicated there are surprising connections between
     alternating permutations and ordinary permutations.

   In this paper, we are concerned with the enumeration of $4123$-avoiding down-up and up-down alternating
   permutations of even and odd length. We  establish    recursive bijections between the set  $\mathcal{DU}_{2n}(4123)$ and
   the set of standard Young tableaux of shape $(n,n,n)$, and between  the set  $\mathcal{DU}_{2n-1}(4123)$ and   the set of shifted  standard
    Young tableaux of shape $(n+1, n, n-1)$ via an intermediate structure of Yamanouchi words. Consequently,  we prove the conjectures, posed by Lewis \cite{Lewis2},
     that $|\mathcal{UD}_{2n}(1432)|=|\mathcal{UD}_{2n}(1234)|$ and $|\mathcal{UD}_{2n+1}(1432)|=|\mathcal{UD}_{2n+1}(2143)|$ in the sense that
     $|\mathcal{UD}_n(1432)|=|\mathcal{DU}_n(4123)|$ by the operation of
     complement.

   Applying the  bijections between $4123$-avoiding  down-up alternating permutations and   standard Young tableaux, we show that
    for $n\geq 1$, $4123$-avoiding  up-down alternating permutations of   $2n+1$  are in one-to-one correspondence with   standard Young tableaux of shape $(n+1,n,n-1)$.
    Moreover,  for $n\geq 2$, $4123$-avoiding up-down alternating permutations of length $2n$  are in bijection  with  shifted  standard Young tableaux of shape $(n+2, n, n-2)$. As a result, we deduce that
 $|\mathcal{UD}_{2n}(4123)|=|\mathcal{UD}_{2n}(1234)| $ and   $|\mathcal{UD}_{2n+1}(4123)|=|\mathcal{UD}_{2n+1}(1234)|$, as conjectured by Lewis \cite{Lewis2}.

 The paper is organized as follows. In Section $2$, we introduce the bijections between the set $\mathcal{DU}_{2n}(4123)$ and the set of standard Young tableaux
 of shape $(n,n,n)$, and between  the set  $\mathcal{DU}_{2n-1}(4123)$ and   the set of shifted  standard Young tableaux of shape $(n+1, n, n-1)$. In section 3,   we get the
  enumeration of $4123$-avoiding up-down alternating permutations of even and odd length.

\section{$4123$-avoiding down-up alternating permutations }
In this section, we aim to construct   recursive bijections between the set  $\mathcal{DU}_{2n}(4123)$ and the set of standard Young tableaux of shape $(n,n,n)$, and between 
 the set  $\mathcal{DU}_{2n-1}(4123)$ and   the set of shifted  standard Young tableaux of shape $(n+1, n, n-1)$.  To do this,   we need
to introduce the following  definitions and notations.

A {\em partition} $\lambda$ of a positive integer $n$ is defined to
be a sequence $(\lambda_1, \lambda_2, \ldots, \lambda_m)$ of
nonnegative  integers such that $\lambda_1+\lambda_2+\ldots
\lambda_m=n$ and $\lambda_1\geq \lambda_2\ldots \geq \lambda_m$.
Given a partition $\lambda=(\lambda_1, \lambda_2, \ldots,
\lambda_m)$, the (ordinary) Young diagram of shape $\lambda$ is the
left-justified array of $\lambda_1+\lambda_2+\ldots+\lambda_m$ boxes
with $\lambda_1$ boxes in the first row, $\lambda_2$ boxes in the
second row, and so on.  If $\lambda$ is a partition with distinct
parts then the {\em shifted} Young diagram of shape $\lambda $ is an
array of cells with $m$ rows, each row indented by one cell to the
right with respect to the previous row, and $\lambda_i $ cells in
row $i$. We shall frequently use the same symbols for things which
may have an ordinary or shifted interpretation. It will always be
clear which interpretation is meant.

If $\lambda$ is a  Young diagram with $n$ boxes, a {\em   standard
Young tableau } of shape $\lambda$ is a filling of the boxes of
$\lambda$ with $[n]$ so that each element appears in exactly one box
and entries increase along rows and columns.   We identify boxes in
Young diagrams and tableaux using matrix coordinates. For example,
the box in the first row and second column is numbered $(1,2)$.

 Given a standard Young
tableau $T$ with $n$ entries, we associate $T$ with a word $\chi(T)$
of length $n$ on the alphabet $\{1,2,3,\ldots\}$, where $\chi(T)$ is
obtained from $T$ by letting the $j$th letter be the row index of
the entry of $T$ containing the number $j$. The words $\chi(T)$ are
known as {\em Yamanouchi} words \cite{eu}. On the other hand, given
a Yamanouchi  word $w$, it is straightforward to recover the
corresponding tableaux $\chi^{-1}(w)$ by letting the $i$-th row
contain the indices of letters of $w$ that are equal to $i$. For
example, the associated standard Young tableau of the Yamanouchi
word $112311223$ is illustrated as follows:
\begin{center}
\parbox[t]{4cm}{ \young(1256,378,49) }
\end{center}

Given a
  word $w$ on the alphabet $\{1,2,\ldots\}$, we define $c_i$ to be the number of   entries
of $w$ that are equal to $i$ and the  {\em type} of the word $w$ to
be the sequence $(c_1, c_2, c_3, \ldots)$. Let $w=w_1w_2\ldots w_n$
be a    word on the alphabet $\{1,2,3\}$. The subsequence
$w_1w_2\ldots w_j$ is said to be an {\em initial run} of length $j$
in $w$ if $w_{j+1}$ is the leftmost entry of $w$ that is equal to
$3$. Similarly,  we can define the {\em final run} of length $j$ to
be a subsequence $w_{n+1-j}w_{n+2-j}\ldots w_{n}$ such that
$w_{n-j}$ is the rightmost entry     equal to $1$.  Denote by
$\alpha(w)$ and $\beta(w)$ the length of the initial run and the
final run of $w$, respectively. For instance, let
$w=121211231323233$. We have $\alpha(w)=7$ and $\beta(w)=6$.

In order to establish the bijections between $4123$-avoiding down-up
 alternating permutations  and   standard Young tableaux, we consider the
following two sets.   Given a  permutation $\pi=\pi_1\pi_2\ldots
\pi_{n}$, let
$$
\mathcal{A}(\pi)=\{0\}\cup \{k|\,\, \exists i<j  \,\,\mbox{ s.t.} \,\,
\pi_i=k, \pi_j={k+1}\, \mbox{and} \,\, k\leq \pi_1-2\}.
$$ Given a   word $w=w_1w_2\ldots w_n$ on the alphabet $\{1,2,3\}$, we define

$$
\mathcal{B}(w)=\{0\}\cup \{k| \,\,w_{k}w_{k+1}=12 \,\,\mbox{and}\,\, k\leq
\alpha(w)-2\}.
$$
\begin{example}
Consider $\pi=658397(10)142$ and $w=121211231323233$. We have
$\mathcal{A}(\pi)=\{0,1,3\}$ and $\mathcal{B}(w)=\{0,1,3\}$.
\end{example}

Given a permutation $\pi\in \mathcal{S}_n$ and an element $a\in
[n+1]$, there is a unique permutation $\pi'=\pi'_1\pi'_2\ldots
\pi'_{n+1}\in \mathcal{S}_{n+1}$ such that $\pi'_1=a$ and the word
$\pi'_2\pi'_3\ldots \pi'_{n+1}$ is order-isomorphic to $\pi$. We
denote this permutation by $a\rightarrow \pi$.  Let $a,b\in [n+2]$
with $b<a$. Denote by $(a,b)\rightarrow \pi$ the permutation
$u=u_1u_2\ldots u_{n+2}$ such that $u_1=a$, $u_2=b$ and
$u_3u_4\ldots u_{n+2}$ is order-isomorphic to $\pi$. More
precisely, the permutation $u$ is defined by
$$
u_i=\left\{\begin{array}{ll} a, & \,\,  i=1,\\
b, & \,\,  i=2,\\
\pi_{i-2}, &\,\, \pi_{i-2}<b,\\
\pi_{i-2}+1, &\,\, b\leq \pi_{i-2}<a-1,\\
\pi_{i-2}+2, &\,\, \pi_{i-2}\geq a-1.
\end{array}\right.
$$

We start  with two lemmas that will be essential in the construction
of the bijections between $4123$-avoiding down-up alternating
permutations and standard Young tableaux. First, we present the
following  simple  observation that will be of use in the subsequent
proofs of Lemmas.

\begin{observation}\label{obser}
Let $\pi=\pi_1\pi_2\ldots \pi_n\in \mathcal{S}_n$ with
  $\mathcal{A}(\pi)=\{a_0, a_1,  \ldots, a_p\}$, where $p\geq 0$ and
$0=a_0<a_1<a_2\ldots<a_p$. Assume that  $a_{p+1}=\pi_1$. For
any integers $r$ and $s$ with $a_j<r<s\leq a_{j+1}$, suppose
that $\pi_l=r$ and $\pi_m=s$. Then we have $l>m$.
\end{observation}

\begin{lemma} \label{lem4123}
Let $\pi=\pi_1\pi_2\ldots \pi_n\in \mathcal{DU}_n(4123)$ with
  $\mathcal{A}(\pi)=\{a_0, a_1,   \ldots, a_p\}$, where $p\geq 0$ and
$0=a_0<a_1<a_2\ldots<a_p$. Assume that  $a_{p+1}=\pi_1$. If
$\pi'=(a,b)\rightarrow \pi $ is a permutation in
$\mathcal{DU}_{n+2}(4123)$, then $b\leq \pi_1$ and there exists an
integer $j$ such that    $a_{j+1}+2\geq a>b\geq a_j+1$.
\end{lemma}
\pf Let $\pi'=\pi'_1\pi'_2\ldots \pi'_{n+2}$.  Recall  that
$\pi'_1=a$ and $\pi'_2=b$. Since $\pi$ is a down-up alternating
permutation, we have $b\leq \pi_1=a_{p+1}$. Suppose that $a_{j+1}\geq b \geq
a_{j}+1$ for some integer $j$ with $0\leq j\leq p$.  We claim that
$a\leq a_{j+1}+2$. Otherwise, assume  that $a> a_{j+1}+2$. Then we
have two cases.   If $j=p$, then the subsequence $ab(a_{p+1}+1)
(a_{p+1}+2)$ is order-isomorphic to  $4123$ in $\pi'$ since $\pi'_3=\pi_1+1=a_{p+1}+1$ and $\pi'_2=b<\pi'_3$. If $j<p$,
 then according to the definition of $\mathcal{A}(\pi)$, there exists
integers $l$ and $m$ with $l<m$ such that $\pi_l=a_{j+1}$ and
$\pi_{m}=a_{j+1}+1$. Note that $\pi'_{l+2}=\pi_l+1=a_{j+1}+1$ and
$\pi'_{m+2}=\pi_m+1=a_{j+1}+2$. Then the subsequence
$\pi'_1\pi'_2\pi'_{l+2}\pi'_{m+2}$ forms a $4123$ pattern in $\pi'$.
Hence, we deduce    that $a_{j}+1\leq b<a\leq a_{j+1}+2$. This
completes the proof. \qed

\begin{lemma}\label{lem2}
Let $\pi=\pi_1\pi_2\ldots \pi_n\in \mathcal{DU}_n(4123)$  with
$\mathcal{A}(\pi)=\{a_0, a_1,  \ldots, a_p\}$, where $p\geq 0$ and $
0=a_0<a_1<a_2\ldots<a_p$. Assume that $a_{p+1}=\pi_1$.  Let
$a,b$ be two integers such that $a_{j+1}+2\geq a>b\geq a_j+1$ and
$b\leq \pi_1$. Then $\pi'=(a,b)\rightarrow \pi$
  is a permutation in
$\mathcal{DU}_{n+2}(4123)$ satisfying that
\begin{itemize}

 \item[{\upshape (\rmnum{1})}]  if $b=a_j+1$ and $j\geq 1$, we have $\mathcal{A}(\pi')=\{a_0, a_1, \ldots, a_{j-1}, b\}$ when $a>b+1$ and $\mathcal{A}(\pi')= \{a_0, a_1, \ldots, a_{j-1}\}$ when $a=b+1$;
\item[{\upshape (\rmnum{2})}]  otherwise, we have $\mathcal{A}(\pi')=\{a_0, a_1, \ldots, a_j, b\}$ when $a>b+1$ and  $\mathcal{A}(\pi')=\{a_0, a_1, \ldots, a_j\}$ when $a=b+1$;
\end{itemize}
\end{lemma}
\pf  Since $b\leq \pi_1$, the permutation $\pi'$ is a down-up alternating permutation.  Now we proceed to show that $\pi'$ avoids the pattern $4123$.
Let $\pi'=\pi'_1\pi'_2\ldots \pi'_{n+2}$. Suppose that there is a
subsequence $\pi'_k\pi'_l\pi'_{m}\pi'_{q}$ with $k<l<m<q$ which is
order-isomorphic to $4123$.
  Since  the subsequence $\pi'_3\ldots \pi'_{n+2}$  is order-isomorphic to $\pi$, we have either $k=1$ or $k=2$.
   If $k=2$, since $\pi'_2<\pi'_3$,   the subsequence  $\pi'_3\pi'_l\pi'_{m}\pi'_{q}$ is  order-isomorphic to the pattern $4123$.
   If $k=1$ and $l> 2$,
     then the subsequence $\pi'_3\pi'_l\pi'_{m}\pi'_{q}$ is an instance of  $4123$ since $\pi'_3\geq \pi_1+1$ and $\pi'_1=a\leq \pi_1+2$.
     Thus, it follows that $k=1$ and $l=2$. Recall that $\pi'_1=a$, $\pi'_2=b$, which  implies that
$ b<\pi'_{m}<\pi'_{q}<a $.  From this, we deduce that $  a_j \leq
b-1<  \pi'_m-1<\pi'_q-1< a-1\leq a_{j+1}+1$.  Note that
$\pi'_{m}=\pi_{m-2}+1$  and $\pi'_{q}=\pi_{q-2}+1$.  So we have
$a_j< \pi_{m-2}<\pi_{q-2}\leq a_{j+1}$. This contradicts with
Observation \ref{obser}.   Thus, the permutation $\pi'$ is in
$\mathcal{DU}_{n+2}(4123)$.

It remains to prove that the permutation $\pi'$\ verifies the points
(\rmnum{1}) and (\rmnum{2}). It is easily seen for any $0\leq k\leq
j-1$, we have $a_k\in \mathcal{A}(\pi')$. Now we proceed to  to show
that there exists no integer $k$ such that $k>b$ and $k\in
\mathcal{A}(\pi')$. Otherwise, suppose that $k$ is such an integer.
According to the definition of $\mathcal{A}(\pi')$, there exists
integers $l$ and $m$ with $l<m$ such that $\pi'_l=k$, $\pi'_m=k+1$
and $k\leq \pi'_1-2=a-2\leq a_{j+1}$. This implies that
$\pi_{l-2}=k-1$ and $\pi_{m-2}=k$. So we have $a_j\leq
b-1<k-1=\pi_{l-2}<\pi_{m-2}=k\leq a_{j+1}$. This contradicts with
Observation \ref{obser}.   So we conclude that there exists no
integer $k$ such that $k>b$ and $k\in \mathcal{A}(\pi')$.

If $a>b+1$, then we have $b\in \mathcal{A}(\pi')$ since $b+1$
appears right to $b$ in $\pi'$ and $b\leq a-2=\pi'_1-2$. If $a=b+1$,
then $b\notin \mathcal{A}(\pi')$ since $b+1$ appears left to $b$ in
$\pi'$. Moreover, when $b=a_{j}+1$ and $j\geq 1$,    we have $a_j\notin
\mathcal{A}(\pi')$ since $a_{j}+1$ appears  left to $a_j$ in $\pi'$.
Otherwise, we have $a_j\in \mathcal{A}(\pi')$.  Hence
(\rmnum{1}) and (\rmnum{2}) are  verified. This completes the proof.
\qed

Now we proceed to construct a  recursive bijection between the set
$\mathcal{DU}_{2n}(4123)$ and the set of Yamanouchi words on the
alphabet $ \{1,2,3\}$ of type $(n,n,n)$.

\begin{theorem}\label{phi}
There is a bijection $\phi$  between  the set
$\mathcal{DU}_{2n}(4123)$ and the set of Yamanouchi words on the
alphabet $ \{1,2,3\}$ of type $(n,n,n)$ satisfying that
$\pi_1=\alpha(\phi(\pi))$ and
$\mathcal{A}(\pi)=\mathcal{B}(\phi(\pi))$  for any  permutation
$\pi=\pi_1\pi_2\ldots \pi_{2n}\in \mathcal{DU}_{2n}(4123)$.
\end{theorem}
\pf Now,  we define a map  $\phi$ from $\mathcal{DU}_{2n}(4123)$ to
the set of Yamanouchi words on the alphabet $\{1,2,3\}$ of type
$(n,n,n)$  in terms of a recursive procedure. For $n=1$, we define
$\phi(21)=123$. It is clear that for $n=1$, the claim holds. Now,
given any permutation $\pi=\pi_1\pi_2\ldots \pi_{2n+2}\in
\mathcal{DU}_{2n+2}(4123)$ for $n\geq 1$, we proceed to construct a
Yamanouchi word $w=\phi(\pi)$. Let $\pi'= \pi'_1\pi'_2\ldots
\pi'_{2n}$ be a $4123$-avoiding down-up alternating permutation and
$\pi=(a,b)\rightarrow\pi'$ where $a=\pi_1$ and $b=\pi_2$. Let
  $\mathcal{A}(\pi')=\{a_0, a_1, \ldots, a_p\}$ with $p\geq 0$ and  $
0=a_0<a_1<a_2<\ldots<a_p$. Assume that  $v=\phi(\pi')=v_1v_2\ldots
v_{3n}$. By the induction hypothesis, $v$ is a Yamanouchi word  on
the alphabet $ \{1,2,3\}$ of type $(n,n,n)$ with the property that
$\alpha(v)=\pi'_1$ and $\mathcal{B}(v)=\mathcal{A}(\pi')$. Assume
that $a_{p+1}=\alpha(v)=\pi'_1$. By Lemma \ref{lem4123},
we have $a_{j+1}+2\geq a>b\geq a_{j}+1$ for some integer $j$ and
$b\leq \pi'_1 =\alpha(v)$.
 Now we proceed to construct a word $w= \phi(\pi)$ on the
alphabet $\{1,2,3\}$ from $v$ by distinguishing the following two
cases.
\begin{itemize}
\item [(\rmnum{1})]
If $a=b+1$, then set $w=v_1v_2\ldots v_{b-1}{\bf 123}v_{b}\ldots
v_{3n}$.
\item[(\rmnum{2})] If $a>b+1$, then set
 $$w=v_1v_2\ldots v_{b-1}{\bf 12}v_{b}\ldots v_{a-2}{\bf 3}v_{a-1}\ldots
v_{3n}.$$
\end{itemize}
Clearly,  the obtained word $w$ is a Yamanouchi word of type
$(n+1,n+1,n+1)$ with an initial run of length $a$, that is,
$\alpha(w)=a$. It remains to show  that
$\mathcal{A}(\pi)=\mathcal{B}(w)$.

For the case $a=b+1$, it is easy to check that when $b=a_j+1$ and $j\geq 1$
$\mathcal{B}(w)=\{a_0, a_1, \ldots, a_{j-1}\}$. Otherwise, we
have $\mathcal{B}(w)=\{a_0, a_1, \ldots, a_{j}\}$.  By Lemma \ref{lem2}, we deduce that
$\mathcal{B}(w)=\mathcal{A}(\pi)$.

For the case $a>b+1$, suppose that $w=w_1w_2\ldots w_{3n+3}$. Since
$w_{b}w_{b+1}=12$ and $b \leq a-2=\alpha(w)-2$, we have $b\in
\mathcal{B}(w)$. Moreover, if $b\geq 2$, we have $b-1\notin \mathcal{B}(w)$ since
$w_b=1$. It remains to show that there exists no integer $k$ with
$k>b$ such that $k\in
  \mathcal{B}(w)$.  Otherwise, assume that there is such an integer $k$.
  According to the definition of $\mathcal{B}(w)$,  we have
  $w_kw_{k+1}=12$ and $k\leq \alpha(w)-2=a-2$. Note that
  $w_{b+1}=2$, which implies that $k\geq b+2$.  So, we have
  $w_{k}w_{k+1}=v_{k-2}v_{k-1}=12$ with $k\leq a-2\leq a_{j+1}$.  This implies that $k-2\in
  \mathcal{B}(v)$. However  we have
   $a_j+1\leq  b\leq k-2\leq a_{j+1}-2$.  This contradicts with the definition of
  $\mathcal{B}(v)$.  Thus we deduce that $\mathcal{B}(w)=\{a_0, a_1, \ldots, a_{j-1}, b\}
 $ when $b=a_j+1$ and $j\geq 1$. Otherwise, we have  $\mathcal{B}(w)=\{a_0, a_1, \ldots, a_j, b\}$.  By Lemma \ref{lem2}, we deduce that
$\mathcal{B}(w)=\mathcal{A}(\pi)$.

We conclude that the the obtained word $w$ is a Yamanouchi word of
type $(n+1,n+1,n+1)$ such that $\alpha(w)=a=\pi_1$ and
$\mathcal{A}(\pi)=\mathcal{B}(w)$.

It is sufficient to construct the inverse mapping of $\phi$ to show
that this is a bijection. Given a Yamanouchi  word $w=w_1w_2\ldots
w_{3n}$ of type $(n, n, n)$ with $\mathcal{B}(w)=\{a_0, a_1,  \ldots,
a_p\}$ with $p\geq 0$ and  $0=a_0<a_1<\ldots<a_p$, we wish to recover a $4123$-avoiding down-up alternating
permutations $\phi^{-1}(w)$ in terms of a recursive procedure.
Assume that $a_{p+1}=\alpha(w)$.  If $w=123$, then
define $\phi^{-1}(w)=21$. Obviously, we have $\alpha(w)=2$ and
$\mathcal{A}(\phi^{-1}(w))=\mathcal{B}(w)=\{0\}$. Clearly, the
claim holds for $n=1$. For $n\geq 2$, set $a=\alpha(w)$.
Now we proceed to associate $w$ with an ordered  pair $(v, b)$ by
the following procedure.
\begin{itemize}
\item[$(a)$]   If $w_{a+2}=3$, then  let $b=a_p$  and     $v$ be a
word obtained from $w$ by removing $w_{b}$, $w_{b+1}$ and $w_{a+1}$
from $w$.

\item[$(b)$] If $w_{a+2}\neq 3$, then find the largest   integer $q$ such that   $q\leq a-1$ and  $w_{q}w_{q+1}=12$.
 Let $b= q$  and $v$ be a
word obtained from $w$ by removing $w_{b}$, $w_{b+1}$ and $w_{a+1}$
from $w$.
 \end{itemize}
Finally, we define $\phi^{-1}(w)=(a,b)\rightarrow \phi^{-1}(v) $.

 Now
we proceed to prove that the map $\phi^{-1}$ is the desired map.
For the case $w_{a+2}=3$, since $w$ is a Yamanouchi word with $w_{a+1}=3$ and $w_{a+2}=3$, there are at least two occurrences of $2$'s left to $w_{a+1}$ and the first  occurrence of $2$ is preceded immediately by an entry $1$.  This guarantees that  there exists at least one subsequence $w_kw_{k+1}=12$ with $k\leq a-2$, that is, $\mathcal{B}(w)\neq \{0\}$. Hence we have $b=a_p>0$.  For the case $w_{a+2}\neq 3$, the property of the Yamanouchi word ensures that there exists at least one subsequence $w_kw_{k+1}=12$ with $k\leq a-1$. Thus, in either case,     the word $v$ is   a
Yamanouchi word of type $(n-1,n-1,n-1)$.

   Suppose that
 $\mathcal{B}(v)=\{c_0, c_1,  \ldots,
c_m\}$ with $m\geq 0$ and  $0=c_0<c_1<\ldots <c_{m}$. Assume that $c_{m+1}=\alpha(v)$. Now we proceed
to show that the obtained permutation is in
$\mathcal{DU}_{2n}(4123)$ satisfying that
 $\mathcal{A}(\phi^{-1}(w))=\mathcal{B}(w)$  and  the first element of
 $\phi^{-1}(w)$  is equal to $a=\alpha(w)$  by considering  the
 following cases.
\begin{itemize}
\item If $w_{a+2}=3$,  then we have $\alpha(v)=a-2$ since $v$ has an initial
run of length $a-2$. In this case, we have $a=\alpha(v)+2=c_{m+1}+2$.
  Moreover, since $b=a_p$,  there is no
subsequence $w_kw_{k+1}=12$ in the subsequence $w_{b+2}\ldots w_{a}$
with $k\leq \alpha(w)-2=a-2$. So, we have $c_m\leq
 b-1$. By the induction  hypothesis,    the permutation $\phi^{-1}(v)$ is
in $\mathcal{DU}_{2n-2}(4123)$ whose first element equals
$\alpha(v)$ and $\mathcal{A}(\phi^{-1}(v))=\mathcal{B}(v)=\{c_0, c_1,  \ldots, c_m\}$.  By Lemma \ref{lem2}, we have $\phi^{-1}(w)=(a,b)\rightarrow \phi^{-1}(v)$ is in $\mathcal{DU}_{2n}(4123)$ since $c_{m+1}+2= a>b\geq c_m+1$ and $b=a_p\leq a-2=\alpha(v)$.     Observe  that $w_{b}=1$. This ensures that $w_{b-1}w_{b}\neq 12$
   and $b-1\notin \mathcal{B}(w)$ when $b\geq 2 $.
 Thus, we derive that if  $c_m=b-1$ and $m\geq 1$ then we have  $\mathcal{B}(v)=\{a_0, a_1, \ldots, a_{p-1}, c_m\}$.
 Otherwise, we have $\mathcal{B}(v)=\{a_0, a_1, \ldots, a_{p-1}\}$. Since $b=a_p\leq \alpha(w)-2=a-2<a-1$, we can verify that $\mathcal{A}(\phi^{-1}(w))=\mathcal{B}(w)$ by Lemma \ref{lem2}.

\item  If $w_{a+2}\neq 3$,  then we have $\alpha(v)\geq  a-1$ since
$v$ has an initial run of length  at least $a-1$. Since $b\leq a-1$, we have $b\leq \alpha(v)=c_{m+1}$. This implies that
there exists an entry $j$ such that $c_{j+1}\geq b\geq c_{j}+1$.
Since there is no subsequence of $w_kw_{k+1}=12$  in the subsequence
$w_{b+2}\ldots w_{a}$, we have $a\leq c_{j+1}+2$.
By the induction  hypothesis,    the permutation $\phi^{-1}(v)$ is
in $\mathcal{DU}_{2n-2}(4123)$ whose first element equals
$\alpha(v)$ and $\mathcal{A}(\phi^{-1}(v))=\mathcal{B}(v)=\{c_0, c_1, \ldots, c_m\}$.  Thus, by Lemma \ref{lem2}, it follows that $\phi^{-1}(w)=(a,b)\rightarrow \phi^{-1}(v)$ is in $\mathcal{DU}_{2n}(4123)$ since $c_{j+1}+2\geq a>b\geq c_j+1$ and $b\leq c_{j+1}\leq c_{m+1}=\alpha(v)$. Note that
$w_{b}=1$. It follows  that $w_{b-1}w_{b}\neq 12$
   and $b-1\neq \mathcal{B}(w)$ when $b\geq 2$.
\begin{itemize}
\item If $b=a_p$,   then for each $  k$ with $0\leq k\leq p-1$,
we have $a_k\in \mathcal{B}(v)$. Thus, we have $\mathcal{B}(v)= \{a_0, a_1,  \ldots,  a_{p-1},  c_j,  c_{j+1}, \ldots,
c_m\}$ when $b=c_j+1$ and $j\geq 1$. Otherwise, we have
$\mathcal{B}(v)=\{a_0, a_1,  \ldots,  a_{p-1},    c_{j+1}, \ldots,
c_m\}$. In this case, since $b=a_p\leq \alpha(w)-2=a-2<a-1$, by  Lemma
\ref{lem2}, we have $\mathcal{A}(\phi^{-1}(w))=\{a_0, \ldots, a_{p-1}, a_p\}=\mathcal{B}(w)$.
\item If $b\neq a_ p$, then for each $  k$ with $0\leq k\leq p$,
we have $a_k\in \mathcal{B}(v)$. Thus, we have $\mathcal{B}(v)=\{a_0, a_1,  \ldots,
a_{p}, c_j,   c_{j+1}, \ldots, c_m\}$ when $b=c_j+1$ and $j\geq 1$. Otherwise,
$\mathcal{B}(v)=\{a_0, a_1,  \ldots,  a_{p},    c_{j+1}, \ldots,
c_m\}$. In this case, since $b\neq a_p$, we have $b=a-1$. By  Lemma
\ref{lem2}, it is easily seen that $\mathcal{A}(\phi^{-1}(w))=\{a_0, a_1, \ldots, a_p\}=\mathcal{B}(w)$.
\end{itemize}
\end{itemize}
    Hence, the map $\phi$
is a desired   bijection. This completes the proof. \qed

\begin{example}
Consider a $4123$-avoiding down-up alternating permutation
$\pi=63758142$. We can obtain a Yamanouchi word $w$ from $\pi$
recursively as follows:
$$
\begin{array}{lllllll}
\pi={\bf 63}758142 &\rightarrow  &{\bf 54}6132 & \rightarrow & {\bf
41}32&
\rightarrow & 21\\
w=12{\bf 12}11{\bf 3}23233 &\leftarrow  & 121{\bf 12}{\bf 3}233 &
\leftarrow &
{\bf 12}12{\bf 3}3& \leftarrow &123,\\
\end{array}
$$
 It is easy to check that the word $w$ has
an initial run of length $6$. Moreover, we have
$\mathcal{A}(\pi)=\{0,1,3\}$ and $\mathcal{B}(w)=\{0,1,3\}$. Conversely,
given a Yamanouchi word $w$, we can recover the $4123$-avoiding
down-up alternating permutation $\pi=63758142$ by reversing the
above procedure.
\end{example}

 For $n\geq 1$, let $w=w_1w_2\ldots w_{3n}$ be a word on the
alphabet $\{1,2,3\}$ of type $(n-1, n, n+1)$. Let $\{a_1, a_2,
\ldots, a_{n-1}\}$, $\{b_1, b_2, \ldots, b_n\}$ and $\{c_1, c_2,
\ldots, c_{n+1}\}$ be the set of indices of letters of $w$ that are
equal to $1$, $2$  and $3$, respectively. Suppose that
$a_1<a_2<\ldots<a_{n-1}$, $b_1<b_2<\ldots <b_{n}$ and
$c_1<c_2<\ldots< c_{n+1}$. If $b_{n}<c_{n}$ and  for all $1\leq j\leq n-1$, we have
$a_j<b_j<c_j$,  then the word $w$ is called a {\em
skew Yamanouchi  word} of type $(n-1,n, n+1)$. For example, let $w=112123231323233$ of type $(4,5,6)$. We have $\{a_1, a_2, a_3, a_4 \}=\{1,2,4,9\}$, $\{b_1, b_2, b_3, b_4, b_5\}=\{3,5,7,11,13\}$ and $\{c_1, c_2, c_3, c_4, c_5, c_6 \}=\{6,8,10,12,14,15\}$. Hence, the word $w$ is a skew Yamanouchi  word of type $(4,5,6)$.

We seek to enumerate $\mathcal{DU}_{2n-1}$ by aping our bijection
$\phi$ for even length permutations. Given a   permutation
$\pi=\pi_1\pi_2\ldots \pi_{2n-1}\in \mathcal{DU}_{2n-1}(4123)$, we
proceed
 to construct a shifted Yamanouchi word  $\psi(\pi)$ on the alphabet $\{1,2,3\}$ of type $(n-1,n,n+1)$. If $n=1$,
 then define $\psi(1)=233$. Since the word $233$ has an initial run of length $1$, we have
 $\alpha(233)=1$. Moreover, we have $\mathcal{A}(1)=\mathcal{B}(233)=\{0\}$.    For $n\geq 2$, set $\psi(\pi)=\phi(\pi)$.   It is easy to check that the arguments in the proof of Theorem \ref{phi} hold for $4123$-avoiding down-up alternating permutation of odd length. As a consequence, we have the following result.

\begin{theorem}\label{psi}
For $n\geq 1$, the map $\psi$ is a bijection between the set
$\mathcal{DU}_{2n-1}(4123)$ and the set of skew Yamanouchi words on
the alphabet $\{1,2,3\}$ of type $(n-1,n,n+1)$ satisfying that that
$\pi_1=\alpha(\psi(\pi))$ and
$\mathcal{A}(\pi)=\mathcal{B}(\psi(\pi))$  for any  permutation
$\pi=\pi_1\pi_2\ldots \pi_{2n-1}\in \mathcal{DU}_{2n-1}(4123)$.
\end{theorem}

For $n\geq 1$, let $w=w_1w_2\ldots w_{3n}$ be a word on the
alphabet $\{1,2,3\}$ of type $(n+1, n, n-1)$.  Let $\{a_1, a_2,
\ldots, a_{n+1}\}$, $\{b_1, b_2, \ldots, b_n\}$ and $\{c_1, c_2,
\ldots, c_{n-1}\}$ be the sets of indices of letters of $w$ that are
equal to $1$, $2$  and $3$, respectively. Suppose that
$a_1<a_2<\ldots<a_{n+1}$, $b_1<b_2<b_3\ldots <b_{n}$ and
$c_1<c_2<\ldots< c_{n-1}$. If $a_2<b_1$ and for all $1\leq j\leq n-1$, we have
$a_{j+2}<b_{j+1}<c_j$, then the word $w$ is called a {\em shifted
Yamanouchi  word} of type $(n+1, n, n-1)$. Note that  for any
shifted standard Young tableau $T$ of shape $(n+1, n, n-1)$, the
word $\chi(T)$ is a  shifted  Yamanouchi  word of type $(n+1, n,
n-1)$. More precisely,  shifted standard Young tableaux  of shape
$(n+1, n, n-1)$ are in bijection with shifted Yamanouchi  words of
type $(n+1, n, n-1)$.

Let $w=w_1w_2\ldots w_{3n}$ be a  word on the alphabet $\{1,2,3\}$. Denote by $rc(w)=(4-w_{n})(4-w_{n-1})\ldots (4-w_{1})$ the
operation of {\em reversed complement}     of $w$.
let $w=w_1w_2\ldots w_{3n}$ be a skew Yamanouchi    word on the
alphabet $\{1,2,3\}$ of type $(n-1, n, n+1)$. Let $\{a_1, a_2,
\ldots, a_{n-1}\}$, $\{b_1, b_2, \ldots, b_n\}$ and $\{c_1, c_2,
\ldots, c_{n+1}\}$ be the set of indices of letters of $w$ that are
equal to $1$, $2$  and $3$, respectively. Obviously, the word $rc(w)$ is a word on the alphabet $\{1,2,3\}$ of type $(n+1, n, n-1)$. For $1\leq i\leq n+1$, set $a'_i=3n+1-c_{n+2-i}$. For  $1\leq i\leq n$, set $b'_i=3n+1-b_{n+1-i}$.  Similarly,  for  $1\leq i\leq n-1$, set $c'_i=3n+1-a_{n-i}$. According to the definition of the reversed complement of $w$, the sets $\{a'_1, a'_2, \ldots, a'_{n+1}\}$, $\{b'_1, b'_2, \ldots, b'_n\}$
and $\{c'_1, c'_2, \ldots, c'_{n-1}\}$ are the sets of indices of letters of $rc(w)$ that are
equal to $1$, $2$  and $3$, respectively. It is easy to check that $rc(w)$ is  a shifted Yamanouchi    word of type $(n+1, n, n-1)$. Indeed,   the operation of    reversed complement   turns out to be a
bijection between   skew Yamanouchi  words of type $(n-1, n, n+1)$
and shifted Yamanouchi  words of type $(n+1, n, n-1)$.   Similarly, the operation of reversed
complement is an involution on the set of Yamanouchi
words of type $(n, n, n)$.  Moreover,
 the operation of reversed complement transforms  an initial run of a word to a final run.  Observe that given  any ordinary or shifted standard Young tableau $T$ of shape $(a,b,c)$ with the $(1,a)$-entry equal to $k$, its corresponding ordinary or shifted Yamanouchi  word
 $\chi(T)$ has a final run of length $a+b+c-k$.
   As an immediate consequence of Theorems \ref{phi} and \ref{psi}, we have the following results.

\begin{theorem}\label{bphi}
The map $\bar{\phi}=\chi^{-1}\circ (rc)\circ\phi$ is a bijection between the set $DU_{2n}(4123)$ and the set of standard Young tableaux of shape $(n,n,n)$ such that the $(1,n)$-entry of the corresponding tableaux is equal to $3n-\pi_1$ for any permutation $\pi=\pi_1\pi_2\ldots \pi_{2n}\in \mathcal{DU}_{2n}(4123)$.
\end{theorem}

\begin{theorem}\label{bpsi}
The map $\bar{\psi}=\chi^{-1}\circ (rc)\circ\psi$ is a bijection between the set $ DU_{2n-1}(4123)$ and the set of shifted  standard Young tableaux of shape
  $(n+1,n,n-1)$  such that   the $(1,n+1)$-entry of the corresponding tableaux is equal to $3n-\pi_1$ for any permutation $\pi=\pi_1\pi_2\ldots \pi_{2n}\in \mathcal{DU}_{2n-1}(4123)$.
\end{theorem}
Recall that there are bijections between the set
$\mathcal{UD}_{2n}(1234)$ and the standard Young tableaux of shape
$(n,n,n)$, and between  the set
     $\mathcal{UD}_{2n+1}(2143)$ and shifted standard Young tableaux of shape $(n+2, n+1,
    n)$. By the operation of complement,  the set
    $\mathcal{DU}_n(4123)$ are in bijection with the set
    $\mathcal{UD}_n(1432)$.  Thus, from Theorems \ref{bphi} and \ref{bpsi},  we derive that $|\mathcal{UD}_{2n}(1432)|=|\mathcal{UD}_{2n}(1234)|$ and
      $|\mathcal{UD}_{2n+1}(1432)|=|\mathcal{UD}_{2n+1}(2143)|$.

\section{4123-avoiding up-down alternating permutations}
In this section, we aim to get the enumeration of $4123$-avoiding
up-down alternating permutations of odd and even length. We will
show that $4123$-avoiding up-down alternating permutations of length $2n+1$ are
in one-to-one correspondence with standard Young tableaux of shape
$(n+1, n, n-1)$. Moreover, for $n\geq 2$, there is a bijection
between the set of  $4123$-avoiding up-down permutations of length
$2n$ and the set of shifted   standard Young tableaux of shape
$(n+2, n, n-2)$. The following Lemma will be essential in
establishing the bijections.
\begin{lemma}\label{4123up}
Let $\pi=\pi_1\pi_2\ldots \pi_n$  be a permutation in
$\mathcal{DU}_n(4123)$ and $a $ be a positive integer. If $a\leq \pi_1$, then
$\pi'=a\rightarrow \pi  $ is in $\mathcal{UD}_{n+1}(4123)$.
\end{lemma}
\pf  Let $\pi'=\pi'_1\pi'_2\ldots \pi'_{n+1}=a\rightarrow \pi$.  In order to prove $\pi'\in \mathcal{UD}_{n+1}(4123)$,
it is sufficient to prove that there exists no subsequence $\pi'_1\pi'_i\pi'_j\pi'_k$ with $i<j<k$ in $\pi'$. Assume  that $\pi'_1\pi'_i\pi'_j\pi'_k$ is a subsequence  order-isomorphic to $4123$. Since $\pi'_1<\pi'_2$,  we deduce that $\pi'_2\pi'_i\pi'_j\pi'_k$ is also a subsequence  order-isomorphic to $4123$, which implies that $\pi_1\pi_{i-1}\pi_{j-1}\pi_{k-1}$ is a subsequence order-isomorphic to $4123$. This contradicts with the fact that $\pi$ is a $4123$-avoiding down-up alternating permutation. This completes the proof. \qed

Now we proceed to construct a map $\gamma$ from the     set
$\mathcal{UD}_{2n+1}(4123)$ to the set of   standard Young tableaux  of shape $(n+1,n,n-1)$. Given a permutation $\pi=\pi_1\pi_2\ldots \pi_{2n+1}\in \mathcal{UD}_{2n+1}(4123)$, define $\pi'=\pi'_1\pi'_2\ldots \pi'_{2n}$ to be a permutation obtained from $\pi$ by removing $\pi_1$ from $\pi$ and deceasing each entry that is larger than $\pi_1$ by one. Obviously,  the permutation $\pi'$ is in $\mathcal{DU}_{2n}(4123)$. By Theorem \ref{bphi}, the tableau  $\bar{\phi}(\pi')$ is a standard Young tableau of shape $(n, n, n)$ with  the $(1,n)$-entry equal to $3n-\pi'_1$.     Define $T=\gamma(\pi)$ to be a tableau obtained from $\bar{\phi}(\pi')$ by deleting the $(3,n)$-entry and inserting a $(1,n+1)$-entry equal to $(3n+1-\pi_1)$. Since $\pi_1\leq \pi_1'$, the obtained tableau $T$ is a standard Young tableau  of shape $(n+1, n, n-1)$. Therefore, the map $\gamma$ is well defined.

\begin{theorem}\label{udodd}
For $n\geq 1$, the map  $\gamma$ is a bijection  between  the set
$\mathcal{UD}_{2n+1}(4123)$ and the set of   standard Young tableaux  of shape $(n+1,n,n-1)$.
\end{theorem}
\pf  It is sufficient to construct the inverse mapping of $\gamma$ to show that $\gamma$ is a bijection. Given a  standard Young tableau  $T$ of shape  $(n+1, n, n-1)$, we wish to recover a permutation $\gamma^{-1}(T)\in \mathcal{UD}_{2n+1}(4123)$. Suppose that the $(1,n+1)$-entry and $(1,n)$-entry of $T$ are equal to  $3n+1-a$ and $3n-b$, respectively. Then we construct a permutation $\gamma^{-1}(T)$ as follows.
\begin{itemize}
 \item Remove the $(1,n+1)$-entry from the tableau $T$ and decrease each entry that is larger than $3n+1-a$ by one;
     \item   Insert a $(3,n)$-entry which is equal to $3n$. Denote by $T'$ the obtained  standard Young tableaux;
\item Finally, set $\gamma^{-1}(T)=a\rightarrow \bar{\phi}^{-1}(T')$.
 \end{itemize}
Note that $T'$ is a standard Young tableau of shape $(n,n,n)$ such that the $(1,n)$-entry equals $3n-b$. Let $\pi'=\bar{\phi}^{-1}(T')=\pi'_1\pi'_2\ldots \pi'_{2n}$.   By Theorem \ref{bphi}, we deduce that $\pi'$ is a down-up alternating permutation in $\mathcal{DU}_{2n}(4123)$ with $\pi_1'=b$.  Since $T$ is a standard Young tableau,  we have $a\leq b$.   By Lemma \ref{4123up}, the obtained permutation $\gamma^{-1}(T)$ is in $\mathcal{UD}_{2n+1}(4123)$.   It is easy to verify that the construction of the map $\gamma^{-1}$ reverses  each step of the construction of the map $\gamma$. This completes the proof. \qed

  Recall that  there is a bijection between the set $\mathcal{UD}_{2n+1}(1234)$  and the set of standard Young tableaux  of shape $(n+1, n, n-1)$ \cite{Lewis1}.
   From Theorem \ref{udodd}, we deduce the following result.
\begin{theorem}
For $n\geq 1$, we have
$$|\mathcal{UD}_{2n+1}(4123)|=|\mathcal{UD}_{2n+1}(1234)|.$$
\end{theorem}

\begin{example}
Consider a $4123$-avoiding up-down alternating permutation
$\pi=4657132$.  Let $\pi'=546132$.  The tableau $\gamma(\pi)$ is illustrated as
\begin{center}
\parbox[t]{4cm}{ \young(1245,369,78) }.
\end{center}

\end{example}

For $n\geq 2$, given a permutation $\pi=\pi_1\pi_2\ldots \pi_{2n}\in \mathcal{UD}_{2n}(4123)$, let $\pi'=\pi'_1\pi'_2\ldots \pi'_{2n}$   be a permutation obtained from $\pi$ by removing $\pi_1$ from $\pi$ and deceasing each entry that is larger than $\pi_1$ by one. Obviously,  the permutation $\pi'$ is in $\mathcal{DU}_{2n-1}(4123)$. By Theorem \ref{bpsi},   the tableau $\bar{\psi}(\pi')$ is a standard Young tableau of shape $(n+1, n, n-1)$ with  the $(1,n+1)$-entry equal to $3n-\pi'_1$.  Finally we obtain a tableau  from $\bar{\psi}(\pi')$ by deleting the $(3,n-1)$-entry and inserting a $(1,n+1)$-entry equal to $(3n+1-\pi_1)$. Since $\pi_1\leq \pi_1'$, the obtained tableau  is a shifted standard Young  tableau  of shape $(n+2, n, n-2)$.
 Therefore we can deduce the following result by the similar arguments as in the proof of Theorem \ref{udodd}.

\begin{theorem}\label{oddud}
For $n\geq 2$,     $4123$-avoiding up-down alternating permutations of length $2n$ are in one-to-one correspondence with  shifted standard Young tableaux of shape $(n+2,n,n-2)$.
\end{theorem}

As in the case for standard Young tableaux, there is a simple hook length formula for shifted standard Young tableaux \cite{ck, stanley1}. By simple computation, we derive that  the number of shifted standard Young tableaux  of shape $(n+2,n,n-2)$ is equal to ${2(3n)!\over n!(n+1)!(n+2)!}.$  Recall that  the number of $1234$-avoiding up-down
    alternating permutations of length $2n$ is given by ${2(3n)!\over
    n!(n+1)!(n+2)!}$.  Hence, we obtain the following result.

\begin{theorem}
For $n\geq 0$, we have $$|\mathcal{UD}_{2n}(4123)|=|\mathcal{UD}_{2n}(1234)|={2(3n)!\over n!(n+1)!(n+2)!}.$$
\end{theorem}

 \noindent{\bf Acknowledgments.}   The   author was supported by
the National Natural Science Foundation of China (no.10901141).


\end{document}